\documentclass{article}
\usepackage{amsmath}
\usepackage{amscd}
\usepackage{amsthm}
\usepackage{amssymb} \usepackage{latexsym}
\usepackage{eufrak}
\usepackage{euscript}
\usepackage{epsfig}
\usepackage{graphics}
\usepackage{array} 
\usepackage{enumerate}
 
\usepackage{pictexwd,dcpic}

\usepackage{color}

\usepackage{hyperref}

\newcommand{\bel}[1]{\begin{equation}\label{#1}}

\newcommand{\be}{\begin{equation}}

\newcommand{\ba}{\begin{eqnarray}}
\newcommand{\ea}{\end{eqnarray}}
\newcommand{\rf}[1]{(\ref{#1})}

\newcommand{\qe}{\end{equation}}
\newcommand{\norm}[1]{\left\lVert#1\right\rVert}

\newcommand{\R}{{\mathbb R}}

\newtheorem{thesis}{Thesis}
\newcommand{\btl}[1]{\begin{thesis}\label{#1}}
\newcommand{\et}{\end{thesis}}

\theoremstyle{theorem}
\newtheorem{theo}{Theorem}[section]
\theoremstyle{corollary}
\newtheorem{coro}{Corollary}[section]
\theoremstyle{lemma}
\newtheorem{lemma}{Lemma}[section]
\theoremstyle{definition}
\newtheorem{defi}{Definition}[section]
\theoremstyle{remark}
\newtheorem*{pf}{Proof}
\theoremstyle{remark}

\theoremstyle{remark}

\title{Topology and curvature of metric spaces}
\author{Parvaneh Joharinad\thanks{p.joharinad@iasbs.ac.ir, Institute for Advanced Studies in Basic Sciences, Zanjan, Iran},   J\"urgen Jost\thanks{jost@mis.mpg.de,  
Max Planck Institute for Mathematics in the Sciences, Leipzig, Germany}}
\begin{document}
\maketitle
\begin{abstract} 
We develop a new concept of non-positive curvature for metric spaces,  based on intersection patterns of closed balls. In contrast to the synthetic approaches of Alexandrov and Busemann, our concept  also applies to metric spaces that might be discrete. The natural comparison spaces that emerge from our discussion are no longer Euclidean spaces, but rather tripod spaces. These tripod spaces include  the
hyperconvex spaces which have trivial  \v{C}ech homology.  This suggests a  link of our  geometrical method to the topological method of  persistent homology employed  in topological data analysis. We also investigate the geometry of general  tripod spaces. 
\end{abstract}


\section*{Introduction}
Curvature is a mathematical notion that originated in classical differential geometry, and through the profound work of Riemann \cite{Riemann:1868}, became the central concept of the geometry of manifolds equipped with a smooth metric tensor, that is, of Riemannian geometry. 
In Riemannian geometry, sectional curvature is a quantity to measure the deviation of a Riemannian metric from being Euclidean, and to compute it, one needs the first and the second derivatives of the metric tensor. However, having an upper or lower bound for the sectional curvature is equivalent to some metric properties which are geometrically meaningful even in geodesic length spaces, that is, in spaces where any two points can be connected by a shortest geodesic. This led to the synthetic theory of curvature bounds in geodesic length spaces or more precisely the Alexandrov \cite{Alexandrov:1957,Berestovskij:Nikolaev:1993} and Busemann \cite{Busemann:1955} definitions of curvature bounds and in particular to the class of spaces with non-positive curvature. These definitions are based on comparing the distance function between two sides and angles in geodesic triangles with their comparison triangles in the Euclidean space $\R^2$. (For the precise definition, see Section \ref{intersection}.) Gromov\cite{Gromov:1999} then systematically explored the concept of non-positive curvature and found various metric characterizations. In particular, the Gromov products will also play a role below.\\

These extensions to metric geometry still require that any two points can be connected by a shortest geodesic. The nonpositive curvature condition can, however, be expressed solely in terms of distances involving quadruples of points. It follows from Reshetnyak's quadrilateral comparison theorem \cite{Reshetnyak:1968} that in a space of nonpositive curvature in the sense of Alexandrov, any quadruple $(x_1,x_2,x_3,x_4)$ of points satisfies the inequality $d^2(x_1,x_4) +d^2(x_2,x_3)-d^2(x_1,x_3)-d^2(x_2,x_4)\le 2d(x_1,x_2)d(x_3,x_4)$ (see for instance Thm.2.3.1 in \cite{Jost:1997}). Berg and Nikolaev \cite{Berg and Nikolaev:2008,Berg and Nikolaev:2018} showed that this condition is in fact equivalent to nonpositive curvature and developed a corresponding theory of nonpositive curvature. 
In a different approach, Alexander et al. \cite{Alexander:2000} used a 4-point criterion that involves an angle comparison in model spaces that applies to all metric spaces and agrees with the Alexandrov comparison in geodesic spaces. This is further explored in \cite{Alexander:2016,Alexander:2019}. In \cite{Jost:et:al: 2015}, a different approach to  curvature inequalities was introduced that works with intersection patterns of distance balls and therefore is also meaningful for discrete metric spaces.  Such intersection patterns have already been investigated from different perspectives. In topological data analysis \cite{Carlsson:2009}, the concept of persistent homology records how such intersection patterns change when the radii of those distance balls increase. This is based on the fundamental idea of \v{C}ech homology which detects nontrivial topology when from a collection of $k$ balls, any $k-1$ of them have a nontrivial  intersection, while the intersection of all $k$ is empty. Such constructions have become useful tools in data analysis, for instance for spike patterns of groups of neurons, see \cite{Curto:2017}. We connect this topological construction here with the notion of a hyperconvex space, that is, a space where any collection of balls with all pairwise intersections being nontrivial also intersects nontrivially. In particular, such a hyperconvex space is topologically trivial in the sense that its \v{C}ech homology is entirely trivial.  \\
Hyperconvex spaces were introduced by Aronszajn and Panitchpakdi \cite{Aronszajn:Panitchpakdi:1956}, who also showed  that every hyperconvex space is complete and  an absolute 1-Lipschitz retract. Hyperconvex spaces can also be characterized as metric spaces with the exact expansion constant equal to 1 (readers are referred  to  section \ref{intersection} for the definition of the expansion constant). Gr\"unbaum studied the relation between the expansion constant  and  the projective constant of normed spaces as well as the relation between the expansion constant and the retraction constant of metric spaces in his PhD thesis. The results of his research are collected in \cite{Grunbaum:1959} without mentioning the connection between hyperconvexity of a metric space and its expansion constant.  Every metric space $X$ has a smallest containing hyperconvex space by a result of Isbell\cite{Isbell:1964}.
Important examples of hyperconvex spaces are metric trees and $L^\infty$-spaces. The former spaces are important for phylogenetic tree reconstruction, as systematically explored by Dress \cite{Dress:1984}. This depends on the observation of  Kuratowski \cite{Kuratowski} that every metric space can be isometrically embedded into an $L^\infty$-space, and his construction will be utilized below. 

Coming back to Riemannian geometry, sectional curvature is a metric invariant of tangent planes. Since in classical geometry, a plane is determined by three points, for a general metric space, sectional curvature should naturally be an invariant for collections of three points. This is the basic motivation for our approach to sectional curvature. Relating it to the discussion of hyperconvexity, we see that we should employ a somewhat weaker property to obtain a model space. We should require that for any collection of three distance balls with the pairwise nonempty intersection, also all three of them should have an intersection point. For geodesic length spaces, this leads us to the notion of a tripod space.
This means that for every triple of points $\{x_1,x_2,x_3\}$ there is a tripod with legs $\{x_1,x_2,x_3\}$. More precisely, see Def. \ref{tripod space} below, a  geodesic length space $(X,d)$ is called a {\em tripod space} if for any three points $x_1,x_2,x_3 \in X$, there exists a point $m\in X$ with
\bel{0a}
d(x_i,m)+d(x_j,m)=d(x_i,x_j),\: \text{for} \: 1\leq i<j\leq 3.
\qe
Hyperconvex spaces are tripod spaces, but since the condition for the latter is weaker than for the former, there are others.

An important point that emerges from our discussion is that the basic model no longer is Euclidean space. Abstractly, curvature expresses the deviation of a given metric space from such a model. Traditional sectional curvature is normalized to be $0$ for Euclidean spaces, and for a Riemannian manifold that is not locally Euclidean it becomes non-zero, but could be either positive or negative. Concerning the generalizations of curvature inequalities due to Alexandrov, Busemann, and others, metric trees have curvature $-\infty$. As discussed, metric trees are special cases of hyperconvex or tripod spaces. In our approach, tripod spaces emerge as the natural models the deviation from which curvature should express. And as argued, this is natural, because such spaces include the hyperconvex ones which have trivial \v{C}ech homology. Therefore, our approach naturally links geometry and topology. \\

Tripod spaces not only include metric trees, but also other classes of spaces that are important in data analysis. 
As a motivating example, consider the  spaces $(\R^n,d_p)$ with $d_p(x,y)=(\sum_{i=1}^n |x^i-y^i|^p)^{1/p}$ for $1\le p <\infty$ and $d_\infty(x,y)=\max_i|x^i-y^i|$. For $1\leq p\leq\infty$ they are geodesic spaces, meaning that any two points $x_0,x_1$ can be connected by a shortest geodesic, that is, a curve whose length is equal to the distance $d_p(x_0,x_1)$. If $1<p<\infty$, the space is uniquely geodesic, that is, such a shortest geodesic is unique, but for $p=1$ and $p=\infty$   that uniqueness of geodesic connections no longer holds for suitable  pairs of points. Take for instance the two points $x=(0,0)$ and $y=(0,1)$ in $(\R^2,d_\infty)$ or $x=(0,0)$ and $z=(1,1)$ in $(\R^2,d_1)$. There are infinitely many  geodesics connecting $x$ and $y$( resp. $x$ and $z$)  with length equal to $1$( resp. 2) with respect to $d_\infty$( resp. $d_1$). The uniqueness of geodesics, however, does not always fail in those spaces. For instance, there is a unique geodesic connecting $x$ to $y$( resp. $x$ to $z$) in the $d_1$( resp. $d_\infty$) metric. The space $(\R^n,d_\infty)$ is a hyperconvex space, hence in particular a tripod space,  and as we will show in this note, $(\R^n,d_1)$ is a tripod space as well, although not hyperconvex for $n>2$ \cite{Jost: 2018}.\\

{\bf Acknowledgement:} We thank Rostislav Matveev for useful discussions. We are also grateful to Stephanie Alexander, I.David Berg  and Urs Lang for useful comments on the first version of our paper.  Parvaneh Joharinad thanks the Max Planck Institute for Mathematics in the Sciences for hospitality and financial support during the work on this paper.  

\section{Intersection of balls}
\label{intersection}
Let $(X,d)$ be a metric space. A continuous mapping $c:[0,1]\longrightarrow X$ is called a {\em path} between $x=c(0)$ and $y=c(1)$ and its length is defined as
\[
\mathit{l}(c):=\sup\sum_{i=0}^{i=n}d(c(t_i),c(t_{i-1})),
\]
where the supremum is taken over the set of all partitions of the interval $[0,1]$.\\
A metric space $(X,d)$ is a {\em length space} if 
\[
d(x,y)=\inf\{\mathit{l}(c):\: \text{c is a path between}\; x \; \text{and}\; y\},
\]
for every $x,y\in X$.\\
A length space $(X,d)$ is said to be {\em geodesic} if each pair of points $x,y\in X$ is connected by a path $c:[0,1]\longrightarrow X$ such that $d(x,y)=\mathit{l}(c)$. Such a path is  called a {\em shortest  geodesic} between $x$ and $y$.\\ 
A point $m\in X$ is called a {\em midpoint} between the points $x,y$ in a metric space $(X,d)$ if $d(x,m)=d(m,y)=\frac{1}{2}d(x,y)$. We say that a pair of points $x,y\in X$ has  {\em approximate midpoints} if for every $\epsilon>0$ there exists $m_\epsilon\in X$ such that
\[
\max\{d(m_\epsilon,x),d(m_\epsilon,y)\}\leq\frac{1}{2}d(x,y)+\epsilon
\]
Every pair of points in a geodesic space (resp. length space) has at least one midpoint (resp. approximate midpoints).\\ 
A metric space $(X,d)$ is {\em convex} if for any distinct points $x_1,x_2\in X$ there is a point $x_0\in X\setminus\{x_1,x_2\}$ with $d(x_1,x_0)+d(x_0,x_2)=d(x_1,x_2)$. Every finitely compact convex metric space is a geodesic space. Convexity can be obtained under the following condition:\\
\textbf{Property A:} For any two points  $x_1,x_2\in X$ and any positive numbers $r_1$ and $r_2$ with $r_1+r_2\geq d(x_1,x_2)$, the intersection $B(x_1,r_1)\cap B(x_2,r_2)$ is nonvoid, where $B(x,r)$ is the closed ball centered at $x$ with radius $r$.\\

\noindent For each pair of points $x_1,x_2$ in the metric space X and any positive numbers $r_1,r_2$ with $r_1+r_2\geq d(x_1,x_2)$, we assign a number 
\[
\rho((x_1,x_2),(r_1,r_2)):=\inf_{x\in X}\max_{i=1,2}\frac{d(x_i,x)}{r_i},
\]
 and define the function $\rho(x_1,x_2)$ as the supremum of these numbers over all positive numbers $r_1,r_2$ with $r_1+r_2\geq d(x_1,x_2)$. If $\rho(x_1,x_2)=1$ for each pair of points $x_1,x_2\in X$, then the existence of approximate midpoints is guaranteed, and  $X$ is a length space provided that it is a complete metric space. If, moreover, the infimum is attained for each pair by some $x_0\in X$, i.e.  \textbf{Property A} is satisfied by $X$, then $X$ is a geodesic space provided that it is complete.\\
To compute $\rho(x_1,x_2)$ in a complete metric space $X$, one only needs to find the infimum  over all $x\in X$ for $r_1=r_2=\frac{1}{2}d(x_1,x_2)$ and hence, 
\be
\rho(x_1,x_2):=\inf_{x\in X}\max_{i=1,2}\frac{d(x_i,x)}{r_i},
\qe
where $r_1=r_2=\frac{1}{2}d(x_1,x_2)$, that is, $\rho(x_1,x_2)=\inf_{x\in X}\max_{i=1,2}\frac{2d(x_i,x)}{d(x_1,x_2)}$. 
\begin{defi}\label{tripod space}
  A geodesic length space $(X,d)$ is a {\em tripod space} if for any three points $x_1,x_2,x_3 \in X$, there exists a median, that is, a point $m\in X$ with
\bel{0}
d(x_i,m)+d(x_j,m)=d(x_i,x_j),\: \text{for} \: 1\leq i<j\leq 3.
\qe
The point $m$ obviously satisfies  $d(x_1,x_2)+d(x_2,x_3)+d(x_3,x_1)=2(d(x_1,m)+d(x_2,m)+d(x_3,m))$.
\end{defi}
\noindent 
Spaces where this median is unique are called median spaces in the literature. For our purposes, however, it is important not to require uniqueness of such a median.\\
The geodesicity assumption can be omitted from the definition, since as we will see in the proof of Theorem \ref{tripod and hyperbolicity}, every tripod space is automatically geodesic.\\ 
The following property which is based on the intersection of three balls guarantees the existence of tripods.\\
\textbf{Property B:} For any three points $x_1,x_2,x_3 \in X$ which are not colinear, i.e. do not lie on a geodesic, and positive numbers $r_1,r_2,r_3$ with $r_i+r_j\geq d(x_i,x_j)$, $1\leq i<j\leq3$, the intersection $\bigcap\limits_{i=1}^{3}B(x_i,r_i)$ is nonvoid.\\

\noindent Similarly, for $x_1,x_2,x_3 \in X$ and positive numbers $r_1,r_2,r_3$ with $r_i+r_j\geq d(x_i,x_j)$ we assign the number
 \[
\rho((x_1,x_2,x_3),(r_1,r_2,r_3)):=\inf_{x\in X}\max_{i=1,2,3}\frac{d(x_i,x)}{r_i},
\]
and define the function $\rho(x_1,x_2,x_3)$ as the supremum of these numbers over all positive numbers $r_1,r_2, r_3$ with $r_i+r_j\geq d(x_i,x_j)$ for $1\leq i<j\leq3$.\\
To compute $\rho(x_1,x_2,x_3)$, one only needs to find the infimum  over all $x\in X$ for $r_1,r_2,r_3$ with $r_i+r_j=d(x_i,x_j)$, $1\leq i<j\leq3$, which is a system of equations with a unique solution obtained by the Gromov products as follows 
\begin{align}\label{Gromov product}
\nonumber r_1&=\frac{1}{2}(d(x_1,x_2)+d(x_1,x_3)-d(x_2,x_3)),&\\
\nonumber r_2&=\frac{1}{2}(d(x_1,x_2)+d(x_2,x_3)-d(x_1,x_3)),&\\
r_3&=\frac{1}{2}(d(x_1,x_3)+d(x_2,x_3)-d(x_1,x_2)).&
\end{align}
Therefore, 
\be \label{rho}
\rho(x_1,x_2,x_3):=\inf_{x\in X}\max_{i=1,2,3}\frac{d(x_i,x)}{r_i},
\qe
where $r_1,r_2,r_3$ are obtained by the relations in \rf{Gromov product}.\\
If $\rho(x_1,x_2,x_3)=1$ and the infimum is attained by some $m\in X$ for the triple $\{x_1,x_2,x_3\}$, then $m$ is the center of a tripod with legs $x_1,x_2,x_3$.\\
 \textbf{Properties A} and \textbf{B} are in fact special cases of the hyperconvexity condition. The hyperconvexity  condition  is concerned with  the intersection of arbitrary families of closed balls and yields the definition of {\em hyperconvex spaces}.
\begin{defi}\label{hyperconvex}
A metric space $(X,d)$ with the following property is called hyperconvex.\\  
\textbf{Property C:} For any family $\{x_i\}_{i\in I}$ of points in $X$ and positive numbers $\{r_i\}_{i\in I}$ with $r_i+r_j\geq d(x_i,x_j)$ for $i,j\in I$, the intersection 
\[
\bigcap\limits_{i\in I}B(x_i,r_i)
\]
is non-empty.
\end{defi}
\noindent We exclude the obvious case where one of the radii is zero. The condition $r_i+r_j\geq d(x_i,x_j)$ can be replaced with the condition $ B(x_i,r_i)\cap B(x_j,r_j)\neq\emptyset$ for all $i,j\in I$, if $X$ is a convex metric space. It is clear that every hyperconvex space is a tripod space.\\
For any family $\{x_i\}_{i\in I}$ of points in $X$ and positive numbers $\{r_i\}_{i\in I}$ with $r_i+r_j\geq d(x_i,x_j)$ for $i,j\in I$, we assign the number
\[
\rho(\{x_i\},\{r_i\}):=\inf_{x\in X}\sup_{i\in I}\frac{d(x_i,x)}{r_i},
\]
and define the function $\rho(\{x_i\})$ as the supremum of these numbers over all positive numbers $\{r_i\}_{i\in I}$ with $r_i+r_j\geq d(x_i,x_j)$ for $i,j\in I$.\\
The least upper bound of the function $\rho$ is equal to the  {\em expansion constant} of the metric space of $(X,d)$ which is defined as the greatest lower bound of numbers $\mu$ for which \[
\bigcap\limits_{i\in I}B(x_i,\mu r_j)\neq\emptyset
\] for any family of closed balls $\{B(x_i,r_i)\}_{i\in I}$ in $X$ with $r_i+r_j\geq d(x_i,x_j)$ for $i,j\in I$.Moreover, if  this infimum is achieved by some $\mu_0$,  then the expansion constant is called exact, c.f. \cite{Grunbaum:1959}.

\section{Hyperconvex spaces}
According  to Definition \ref{hyperconvex}, a geodesic metric space $(X,d)$ is hyperconvex if every family of closed balls $\{B(x_i,r_i)\}_{i\in I}$ in $X$ with non-empty pairwise intersections, has a non-empty intersection itself. In other words, for a hyperconvex space, the expansion constant is $1$ and is exact.\\
Hyperconvex spaces have the following properties
\begin{enumerate}
\item[a)] Hyperconvex spaces  are complete and  contractible to each of their points 
\cite{Aronszajn:Panitchpakdi:1956, Isbell:1964}.
\item[b)] The metric space $X$ is hyperconvex if and only if every $1-$Lipschitz map from a subspace of any metric space $Y$ to $X$ can be extended to a $1-$Lipschitz map over $Y$ \cite{Aronszajn:Panitchpakdi:1956}.
\item[c)] Isbell \cite{Isbell:1964} proved that every metric space is isometrically embedded in a hyperconvex space, which is unique up to isometry, called its injective envelope or injective hull. The injective envelope of a compact space is compact and the injective envelope of a finite space is a simplicial complex \cite{Isbell:1964}. This construction  was rediscovered by Dress  \cite{Dress:1984}), who called it the  tight span, and who further investigated its  combinatorial properties.  In the hands of Dress and his coworkers, this construction became an important tool in the reconstruction of phylogenetic relationships, see \cite{Dress et al:2002,Dress et al:1996}.  Because of the results of Isbell and Dress, the injective envelope is also called the hyperconvex hull in the literature, and we shall also use that terminology below. 
 \item[e)] Any hyperconvex space included in a metric space is a non-expansive retraction of this space and hyperconvexity is preserved under non-expansive retractions, i.e. the metric space $X'$ is hyperconvex if there is a non-expansive contract from the hyperconvex space $X$ to $X'$ \cite{Grunbaum:1959}.
\item[f)] Any family of hyperconvex spaces with non-empty finite intersection has non-empty intersection.
\item[g)] A complete metric space $Y$ is hyperconvex if for each $\epsilon>0$ there exists a hyperconvex subspace $X$ of $Y$ such that $Y$ is contained in the $\epsilon-$neighbourhood of $X$.
\end{enumerate}
Let $\{(X_i,d_i)\}_{i\in I}$ be a family of hyperconvex spaces and for each $i\in I$ fix a point $z_i\in X_i$. The $\mathit{l}_\infty$ product of this family with the chosen fixed points, i.e. the set of all $(x_i)_{i\in I}$ with $x_i\in X_i$ and $\sup_{i\in I}d_i(x_i,z_i)<\infty$ endowed with the metric $d_{\infty}((x_i),(y_i)):=\sup_{i\in I}d_i(x_i,y_i)$, is a hyperconvex space. If for an arbitrary index set $I$ one takes $\R$ as $X_i$ and $0$ as the fixed point for all $i\in I$, then the corresponding $\mathit{l}_{\infty}$ product is the Banach space $\mathit{l}_{\infty}(I)$ which is then a hyperconvex space since $\R$ is hyperconvex according to the item b above.  $\mathit{l}_\infty(I)$ is the set of all bounded functions  $f:I\longrightarrow \R$ with metric given by the  sup norm, i.e. $d_\infty(f,g):=\sup_{i\in I}|f(i)-g(i)|$. For instance, when $X$ is a finite set consisting of $n$ elements, then $\mathit{l}_{\infty}(X) = (\R^n, d_{\infty})$.\\
A natural way to isometrically embed the metric space $(X,d)$ into a hyperconvex space is by the canonical embedding  $x\mapsto d_x-d_z$ of $X$ in $l_{\infty}(X)$, where $d_x(y):=d(x,y)$ and $z\in X$ is any pre-chosen base point. $X$ is  a hyperconvex space itself  if the image of this embedding is a $1$-Lipschitz retract in $l_{\infty}(X)$, according to the item e above.\\

Let $(X,d)$ be a geodesic space. A geodesic triangle on $X$ consists of three points $x_1,x_2,x_3\in X$ and three minimal geodesics $c_{ij}:[0,1]\longrightarrow X$, $1\leq i<j\leq 3$, parametrized proportional to arc length with $c_{ij}(0)=x_i$ and $c_{ij}(1)=x_j$. We will refer to such a triangle as $\bigtriangleup(x_1,x_2,x_3)$ without specifying its geodesic sides, although they need not be unique and the triangle will therefore depend on the choice of these geodesics. 
\begin{defi}\label{delhyp}
The geodesic space (X,d) is $\delta$-hyperbolic if every geodesic triangle is $\delta$-slim, which means that each side of every geodesic triangle is in the $\delta$-neighbourhood of the union of the two other sides.
\end{defi}
There is another definition of $\delta$-hyperbolicity (also called Gromov-$\delta$-hyperbolicity) that is often employed in the literature. 
\begin{defi}\label{delhypalt} 
The geodesic space (X,d) is $\delta$-hyperbolic if for any collection $x_1,x_2,x_3,x_4\in X$
\begin{equation*}
  d(x_1,x_2)+d(x_3,x_4)\le \max \{d(x_1,x_3)+d(x_2,x_4),d(x_1,x_4)+d(x_2,x_3)\} +\delta.
\end{equation*}
\end{defi}  
The Definitions \ref{delhyp} and \ref{delhypalt} are essentially equivalent.  In fact, If $X$ is $\delta$-hyperbolic according to Def. \ref{delhyp} (triangles are $\delta$-slim) then it is $8\delta$-hyperbolic according to Def. \ref{delhypalt}  (Gromov $8\delta$-hyperbolic). Conversely, if $X$ is Gromov $\delta$-hyperbolic, then it has $4\delta$-slim triangles.\\
Here, however, we shall work with Def. \ref{delhyp}.\\
For each geodesic triangle $\bigtriangleup(x_1,x_2,x_3)$ in the geodesic space $(X,d)$ with sides $c_{12}$, $c_{13}$ and $c_{23}$, there are three points $m_{12}$ on $c_{12}$, $m_{13}$ on $c_{13}$ and $m_{23}$ on $c_{23}$ with 
\begin{equation*}
d(x_i,m_{ij})=r_i \: \text{and}\: d(x_j,m_{ij})=r_j, \; \text{for} \; 1\leq i< j\leq 3, 
\end{equation*} 
where $r_1,r_2,r_3$ are the Gromov products defined in  \rf{Gromov product}. The {\it insize} of $\bigtriangleup(x_1,x_2,x_3)$ is defined as the diameter of the set $\{m_{12},m_{13},m_{23}\}$ and $(X,d)$ is $\delta$-hyperbolic provided that the insize of each geodesic triangle on $X$ is less than or equal to $\delta$. \\  
The set of hyperconvex spaces includes complete $\R$-trees, the only $0$-hyperbolic geodesic spaces, but there are other hyperconvex spaces besides $\R-$trees. For instance,  the space $(\R^n,d_\infty)$ is  a hyperconvex space which is not $\delta$-hyperbolic for any $\delta$. However, A. Dress \cite{Dress:1984} (see also \cite{Dress et al:1996}) showed that the hyperconvex hull (tight extension in his terminology)  of a $0$-hyperbolic space in the sense of Def.\ref{delhypalt}  is geodesic and $0$-hyperbolic. The corresponding for  $\delta$-hyperbolic spaces was stated without proof in \cite{Dress et al:1996}, as the proof in \cite{Dress:1984} can be extended. A proof of the general result with a different method was provided by U. Lang in \cite{Lang:2013} where it is  also proved that the hyperconvex hull of a geodesic $\delta$-hyperbolic space is within distance $\delta$  from that space.

Apparently unaware of Dress' work, Bonk and Schramm \cite{Bonk and Schramm:2000} gave a constructive proof that any $\delta$-hyperbolic
metric space X embeds isometrically in a complete geodesic $\delta$--hyperbolic
metric space.
 The following corollary can also be used to show that every geodesic $0$-hyperbolic space is hyperconvex.
 \begin{coro}
If $(X,d)$ is a geodesic $\delta$-hyperbolic space and  $\{B(x_i,r_i)\}_{i\in I}$  a family of closed balls in $X$ with $r_i+r_j\geq d(x_i,x_j)$, then 
\bel{del1}
\bigcap_{i\in I}B(x_i,\delta+r_i)\neq\emptyset.
\qe 
\end{coro}
 \begin{pf}
Let $E(X)$ be the injective envelope of X. We already know that $E(X)$ is hyperconvex and at distance $\delta$ from $X$. Therefore, the family $\{B(x_i,r_i)\}_{i\in I}$ has a common point in $E(X)$, say $m$, which is at distance at most $\delta$ from $X$. So there exists $x\in X$ such that $d(x,m)\leq\delta$ and consequently $x$ belongs to the intersection $\bigcap_{i\in I}B(x_i,\delta+r_i)$ by the triangle inequality. \qed
 \end{pf}

Hyperconvex spaces are special cases of more general spaces called $\lambda$-hyperconvex spaces.
\begin{defi}
A metric space $(X,d)$ is said to be $\lambda$-hyperconvex for $\lambda\geq1$ if for every family $\{B(x_{i}, r_{i})\}_{i\in I}$ of closed balls in $X$ with the property $r_{i}+ r_{j}\geq d(x_{i}, x_{j})$, one has
\bel{del2}\bigcap\limits_{i\in I}B(x_{i},\lambda r_{i})\neq\emptyset.
\qe
\end{defi}
The difference between \rf{del1} and \rf{del2} makes the difference between the corresponding concepts clear. For large radii, the $\delta$ in \rf{del1} becomes insignificant, and so, this is good for asymptotic consideration, while the condition \rf{del2} is invariant under scaling the metric $d$.\\
\noindent It is also clear that if the expansion constant of the metric space $X$ is equal to $\mu$, then X is $\lambda$-hyperconvex for every $\lambda>\mu$ and $X$ is hyperconvex if and only if it is $1$-hyperconvex.\\
Hilbert spaces are $\sqrt{2}$-hyperconvex. Reflexive Banach spaces and dual Banach spaces are $2$-hyperconvex. Therefore, for the measure space $(X,\mu)$, the spaces $L^p(X,\mu), 1<p<\infty$, are $2$-hyperconvex, and if $X$ is finite, the  finite dimensional space $L^1(X,\mu)$  is also $2$-hyperconvex. $L^\infty(X,\mu)$ is $1$-hyperconvex, as pointed out above.

\section{Upper curvature bounds}\label{upper}

There are different formulations of non-positively curved geodesic spaces, the best known ones being those of  Alexandrov and Busemann. 
\\
After recalling these definitions (see \cite{Jost:1997} for a systematic presentation), we will introduce a new definition of non-positively curved metric spaces, which is based on the intersection property of closed balls. We will show that $CAT(0)$ spaces, i.e. the non-positively curved geodesic spaces in the sense of Alexandrov, have non-positive curvature according to this definition.\\
Our definition is more general because it is meaningful for  general metric spaces. Moreover,
the non-positively curved geodesic spaces in the sense of this new formulation are no longer necessarily uniquely geodesic. For instance, there are hyperconvex spaces or more generally tripod spaces in the category of non-positively curved spaces which are not uniquely geodesic, like the $L^\infty$-spaces already discussed. Therefore,  non-positively curved spaces in our sense need not be non-positively curved in the sense of Busemann, and, a fortiori, not in the sense of Alexandrov either.\\  
Let $(X,d)$ be a geodesic space and $\bigtriangleup(x_1,x_2,x_3)$ a geodesic triangle in $X$.
Then a comparison triangle for $\bigtriangleup(x_1,x_2,x_3)$ is a triangle with vertices $\{\bar{x}_1,\bar{x}_2,\bar{x}_3\}$ in $\R^2$  that satisfies
\be\nonumber
\norm{\bar{x}_i-\bar{x}_j}=d(x_i,x_j),\: 1\leq i<j\leq 3,
\qe  
where $\norm{.}$ denotes the Euclidean norm on $\R^2$.
 \begin{defi}
The geodesic space $(X,d)$ is a $CAT(0)$ space if each pair of geodesics $c_1, c_2:[0,1]\longrightarrow X$ emanating from a point satisfies the following inequality
 \begin{equation}\label{CAT(0)}
 d(c_1(t),c_2(s))\leq \norm{\bar{c}_1(t)-\bar{c}_2(s)},\, \forall\ t,s\in[0,1],
 \end{equation}
where $\bar{c}_1,\bar{c}_2:[0,1]\longrightarrow \mathbb{R}^2$ are the respective sides of the comparison triangle in $\mathbb{R}^2$ for $\bigtriangleup(c_1(0),c_1(1),c_2(1))$.\cite{Burago:Burago:Ivanov:2001}
\end{defi}
\noindent Loosely speaking, triangles in $CAT(0)$ spaces are not thicker than Euclidean triangles with the same side lengths.  A simply connected complete $CAT(0)$ space is called a Hadamard space.\\
In the same manner, one can define $CAT(\kappa)$ spaces by comparing the well defined triangles in these spaces with the corresponding comparison triangles in the model surface $M_{\kappa}^{2}$, i.e. the simply connected complete $2$-dimensional Riemannian manifold with  constant sectional curvature $\kappa$.\\
Let $(X,d)$ be a metric space, $z,x_1,x_2\in X$ such that $x_1,x_2\in X\setminus \{z\}$, and let $\bigtriangleup(\bar{z},\bar{x}_1,\bar{x}_2)$ be the comparison triangle in $\mathbb{R}^2$ for $\bigtriangleup(z,x_1,x_2)$. The comparison angle at $z$ denoted by $\bar{\angle}_{z}(x_1,x_2)$ is by definition the Euclidean angle $\angle_{\bar{z}}(\bar{x_1},\bar{x_2})$. If $X$ is a $CAT(0)$ space and $c_1$ and $c_2$ are geodesics from $z$ to $x_1$ and $x_2$ respectively, the inequality \rf{CAT(0)} implies that
\begin{equation}\label{angle comparison}
 \bar{\angle}_{z}(c_1(t),c_2(t'))\leq\bar{\angle}_{z}(c_1(s),c_2(s')),\, \forall\ 0<t\leq s\leq 1, 0<t'\leq s'\leq 1.
 \end{equation}
This helps to define the angle $\angle_{z}(x_1,x_2)$ by the following limit
\[\angle_{z}(x_1,x_2):=\lim_{t,t'\rightarrow 0}\bar{\angle}_{z}(c_1(t),c_2(t'))\]
The relation \rf{angle comparison} implies that $\angle_{z}(x_1,x_2)\leq \bar{\angle}_{z}(x_1,x_2)$ and therefore the Euclidean law of cosines implies the following cosine inequality in complete $CAT(0)$ spaces
\bel{cosine}
d(x_{1},x_{2})^2\geq d(x_{1},z)^2+d(z,x_{2})^2-2d(x_{1},z)d(z,x_{2})\cos (\angle_{z}(x_{1},x_{2})).
\qe
Busemann's definition of non- positively curved spaces  is  more general. 
 \begin{defi}
 A geodesic space $(X,d)$ is a Busemann convex space (or of non-positive curvature in the sense of Busemann) if for every two geodesics $c_1, c_2:[0,1]\longrightarrow X$ with $c_1(0)=c_2(0)$, the distance function $t\mapsto d(c_1(t),c_2(t))$ is convex. 
 \end{defi}
This property in particular implies that Busemann convex spaces are uniquely geodesic. \\
Every $CAT(0)$ space is Busemann convex but not conversely.\\

We shall now introduce our definition of non-positively curved metric spaces and state  an axiomatic curvature inequality which is  applicable to arbitrary metric spaces. Our definition is based on the intersection of three balls and tripod spaces are the extreme spaces to which other spaces will be compared. As a motivation for considering triples of points, we note  that, in Riemannian geometry, sectional curvatures are associated to tangent planes and  a plane is determined by three points. This therefore  suggests an abstract approach towards
the curvature of general metric spaces by looking at the metric relations between three
points. \\
\noindent We recall from \rf{rho} that for each triple of points $(x_1,x_2,x_3)$ in the metric space $(X,d)$, 
 \be
\rho(x_1,x_2,x_3):=\inf_{x\in X}\max_{i=1,2,3}\frac{d(x_i,x)}{r_i},
\qe
where $r_1,r_2,r_3$ are obtained by the relations in  \rf{Gromov product}.\\
 \begin{defi}\label{proposed definition}
The metric space $(X,d)$ is said to be of non-positive curvature if for each triple $(x_1,x_2,x_3)$ in $X$ with the comparison triangle  $\bigtriangleup(\bar{x}_1,\bar{x}_2,\bar{x}_3)$ in $\R^2$, one has
\be
\rho(x_1,x_2,x_3)\leq \rho(\bar{x}_1,\bar{x}_2,\bar{x}_3),
\qe
where $\rho(\bar{x}_1,\bar{x}_2,\bar{x}_3)$ is similarly defined by 
\be\nonumber
\rho(\bar{x}_1,\bar{x}_2,\bar{x}_3):=\min_{x\in \R^2}\max_{i=1,2,3}\dfrac{\norm{x-\bar{x}_i}}{r_i}.
\qe
\end{defi}
\noindent This condition embodies the principle that for any triple of closed balls $\{B(x_{i}, r_{i});\:i=1,2,3\}$ with pairwise intersection in the geodesic metric space $X$, the intersection $\bigcap_{i=1,2,3}B(x_{i},\rho r_{i})$ is non-empty whenever $B(\bar{x}_{i},\rho r_{i})$, $i=1,2,3$, have a common point.\\
\noindent For every triangle $(x_1,x_2,x_3)$ in the metric space $X$ that is non-positively curved in the sense of  Definition \ref{proposed definition}, the quantity $\rho(\bar{x}_1,\bar{x}_2,\bar{x}_3)$ is smaller than or equal to $\frac{2}{\sqrt{3}}$, i.e. the expansion constant of $\R^2$. By  Helly's Theorem, the intersection of a family of closed balls in $\R^2$ is non-empty if the intersection of every triple of balls in that family is non-empty.\\

We observe
\begin{lemma}\label{tripodlemma}
  A tripod space has non-positive curvature in the sense of  Definition \ref{proposed definition}).
\end{lemma}
\begin{pf}
 By taking the point $m$ to be the center of a tripod with legs $x_1,x_2$ and $x_3$ as  in Definition \ref{tripod space}, we see that 
\bel{tri1}
\rho(x_1,x_2,x_3)=1
\qe
for any three points in a tripod space. \qed
\end{pf}
The following Lemma is a variant of the Lemma 4.2 in \cite{Lang:Schroeder:1997} which is suitable for the current discussion.
\begin{lemma}\label{Hadamard}
Let $(X,d)$ be a complete $CAT(0)$ space, $F=\{x_{1},...,x_{n}\}$ a finite subset of $X$ and $r_{1},..., r_{n}$ positive numbers. Then $\mu':=\inf\{\mu>0: \bigcap_{i=1}^{n}B(x_{i},\mu r_{i})\neq\emptyset\}$ is finite and $\bigcap_{i=1}^{n}B(x_{i},\mu' r_{i})=\{m\}$ for some $m\in\overline{G(F')}$, where $F':=\{x_{i}\in F: d(m,x_i)=\mu' r_i\}$ and $G(F')$ is the convex hull of $F'$.
\end{lemma}
\noindent We present a shorter version of the proof under the extra condition that $X$ is simply connected.
\begin{pf}
Let $(X,d)$ be a Hadamard space  and  $F=\{x_{1},...,x_{n}\}$ a finite subset of $X$. The function  $d^2(x_i,.)$ is strongly convex or a $\lambda$-convex function with $\lambda=1$ for each $1\leq i\leq n$. More precisely, for each unit speed geodesic $c$, the function $t\mapsto d^2(x_i,c(t))-\lambda t^2$ is convex or equivalently  $d^2(x_i,c(t))$ is no less convex than the quadratic function $t\mapsto\lambda  t^2$. The proof is just a reformulation of the definition of non-positive curvature in the sense of Alexandrov and can be found in \cite{Burago:Burago:Ivanov:2001}, Theorem 9.2.19. 
Hence the function 
\begin{equation*}
\max _{i=1,2,...,n}\frac{d^2(x_i,.)}{r_i}
\end{equation*}
has a unique minimizer $m$, since a  continuous proper $\lambda$-convex function on a complete metric space has a unique minimum point.\\
Furthermore,  let $m'\in \overline{G(F')}$ be the closest point to $m$. The  existence and uniqueness of $m'$ is also guaranteed by the strong convexity of the distance function. Let $c:[0,d(m,m')]\longrightarrow X$ be the unique geodesic connecting $m$ to $m'$.  It follows by the choice of $m'$ that $d(c(s), x_i)\leq \mu' r_i$ for all $s\in[0,d(m,m')]$ and $x_i\in F'$, where  
\begin{equation*}
\mu'=\max _{i=1,2,...,n}\frac{d^2(x_i,m)}{r_i}.
\end{equation*} 
Since $d(m,x_j)<\mu'r_j$ for all $x_j\in F\setminus F'$, points on $c$ which are close enough to $m$ coincide with $m$ by the uniqueness of $m$ and therefore $m=m'$. \qed
\end{pf}
\begin{coro}\label{cor1}
Let $(X,d)$ be a complete $CAT(0)$ space, $\{x_1,x_2,x_3\}$ an arbitrary triple of points in $X$ and $r_1,r_2,r_3$ the respective Gromov products as defined in \rf{Gromov product}. Then the infimum in \rf{rho}, which we denote by $\rho$, is attained by a unique point $m$. Moreover,  $d(x_i,m)=\rho r_i$ for $i=1,2,3$.
\end{coro}
\begin{pf}
Without loss of generality, we can suppose that $F'=\{x_1, x_2\}$, since $F'$ has at least two points. Then $m$ lies on the geodesic connecting $x_1$ and $x_2$ and  $d(x_i,m)=\rho r_i$ for $i=1,2$, where $r_1+r_2=d(x_1,x_2)$. Therefore, $\rho=1$ and $m$ is the center of a tripod with legs $x_1,x_2$ and $x_3$, which contradicts the assumption that $x_3$ does not belong to $F'$.\qed
\end{pf}

\begin{theo}\label{Theorem}
Every complete $CAT(0)$ space is $\sqrt{2}$-hyperconvex.
\end{theo}
This result can be derived from Lemma 3.1 in \cite{Lang:2000}. For completeness, we provide the simple proof.
\begin{pf}
Let $\{B(x_{i}, r_{i})\}_{i\in I}$ be a family of closed balls in the complete $CAT(0)$ space $X$ with the non-empty pairwise intersection, which implies that $ r_{i}+ r_{j}\geq d(x_{i}, x_{j})$ for $i,j\in I$. Since $CAT(0)$ spaces have the finite intersection property, it suffices to prove that the intersection of any finite number of balls in $\{B(x_{i},\sqrt{2} r_{i})\}_{i\in I}$ is non empty, i.e. $\bigcap_{i=1}^{n}B(x_{i},\sqrt{2} r_{_i})\neq\emptyset$. Let
$$\mu':=\inf\{\mu>0: \bigcap_{i=1}^{n}B(x_{i},\mu r_{i})\neq\emptyset\}$$
By Lemma \ref{Hadamard},  $\mu'<\infty$ and $\bigcap_{i=1}^{n}B(x_{i},\mu' r_{i}) $ consists of a single point $z$. After a possible reordering, one can assume there is a number $2\leq k\leq n$, for which
\[
\begin{cases}
d(x_{i}, z)=\mu' r_{i}, \, i=1,...,k,&\\
d(x_{i}, z)<\mu' r_{i},\, i=k+1,...,n.&
\end{cases} \]
Moreover, there exists at least one $2\leq j\leq k$ such that $\angle_{z}(x_{1},x_{j})\geq\frac{\pi}{2}$. Assume  $\angle_{z}(x_{1},x_{j})<\frac{\pi}{2}$ for all $2\leq j\leq k$ and $c_{i}:[0,\mu'r_{i}]\longrightarrow X$ are unit speed geodesics from $z$ to $x_{i}$ for $1\leq i\leq k$.  One can choose $t$ small enough such that $d(c_{1}(t), x_{i})<d(z,x_{i})=\mu'r_i$ for $i=1,..,k$ and $c_1(t)$ still lies in the ball $B(x_i,\mu'r_i)$ for $i=k+1,...,n$, which contradicts the uniqueness of $z$.\\
Assuming $\mu'>\sqrt{2}$ and using the cosine inequality \rf{cosine} for $CAT(0)$ spaces, one has
\begin{align*}
d(x_{1},x_{j})^2&\geq d(x_{1},z)^2+d(z,x_{j})^2-2d(x_{1},z)d(z,x_{j})\cos (\angle_{z}(x_{1},x_{j}))&\\
&\geq d(x_{1},z)^2+d(z,x_{j})^2&\\
&=\mu'^2(r_{1}^2+r_{j}^2)&\\
&>2((r_{1}^2+r_{j}^2)&\\
&\geq(r_{1}+r_{j})^2.&
\end{align*} 
This contradicts the fact that $d(x_{1},x_{j})\leq r_{1}+r_{j}$.\qed
\end{pf}
\begin{theo}\label{catgivesnpc}
Every complete $CAT(0)$ space has non-positive curvature in the sense of Definition \ref{proposed definition}
\end{theo}
\begin{pf}
Let $(X,d)$ be a complete $CAT(0)$ space, $\{x_1,x_2,x_3\}$ an arbitrary triple of points in $X$ and $r_1,r_2,r_3$ the respective Gromov products as defined in \rf{Gromov product}. Then according to Lemma \ref{Hadamard} and its  Corollary \ref{cor1}, the infimum in \rf{rho}, which we denote by $\rho$, is attained by a unique point $m$ and $d(x_i,m)=\rho r_i$ for $i=1,2,3$.\\
Let $\bigtriangleup(\bar{x_1},\bar{x}_2,\bar{x}_3)$ be the comparison triangle in $\R^2$ for $\bigtriangleup(x_1,x_2,x_3)$ and $f:\{\bar{x_1},\bar{x}_2,\bar{x}_3\}\longrightarrow X$ be the isometric map with $f(\bar{x}_i)=x_i$ for $i=1,2,3$. Then there exists a $1$-Lipschitz extension $\bar{f}:\R^2\longrightarrow X$ of $f$ according to Theorem A in \cite{Lang:Schroeder:1997}. For the point $\bar{m}=\arg\min\max_{i=1,2,3}\dfrac{\norm{x-\bar{x}_i}}{r_i}$ we have 
\be\nonumber
d(\bar{f}(\bar{m}),x_i)\leq \norm{\bar{m}-\bar{x}_i}=\bar{\rho}r_i,\: i=1,2,3, 
\qe
where $\bar{\rho}=\rho(\bar{x}_1,\bar{x}_2,\bar{x}_3)$. This means that $\bigcap\limits_{i=1,2,3}B(x_i,\bar{\rho}r_i)$ is nonvoid and consequently $\rho\leq\bar{\rho}$.\qed
\end{pf}
The converse of Theorem \ref{catgivesnpc} does not hold, because in non-positively curved spaces in the sense of Definition \ref{proposed definition}, like $(\R^n,d_\infty)$, shortest geodesics need not be unique, while they are in $CAT(0)$-spaces (see e.g. \cite{Berestovskij:Nikolaev:1993,Jost:1997}). When we assume that the metric space in question is a Riemannian manifold, however, we do have the converse statement.
\begin{theo}\label{npcgivessec}
If the complete Riemannian manifold $(N,g)$ has non-positive curvature in the sense of Definition \ref{proposed definition}, then it is of non-positive sectional curvature.
\end{theo}
\begin{pf}
For a triangle with vertices $x_1,x_2,x_3$  in $N$, we denote by $m$ the minimizing point in the definition of $\rho(x_1,x_2,x_3)$, which  is unique when the three points are sufficiently close, the situation we shall only need to consider. \\
Let $V_1,V_2,V_3$  be the tangent vectors of the three shortest geodesics from $m$ to the corners of the triangle at $m$.
 It follows from the minimizing property that $V_1,V_2,V_3$  lie in a plane. In fact, if they did not lie in a plane in $T_mN$, the sum of the three angles at $m$ between them would be smaller than $2\pi$. We could then find a vector $V\in T_mN$ with $\langle V_j,V\rangle >0$ for $j=1,2,3$. Moving $m$ in the direction of $V$,  decreases all three  distances $d(m,x_j)$, hence a contradiction to the definition of  $m$.\\ 
For a close enough triple $(a,b,c)\in N$, let $m$ be the unique point which gives the minimum value $\rho:=\rho(a,b,c)$  and let $V_a,V_b,V_c$ be the tangent vectors at $m$ to the shortest geodesics from $m$ to the vertices $a,b,c$ respectively. Take a triple $(A,B,C)\in \R^2$ with $\|A\| =d(a,m)$ etc,  and the same angles between the shortest geodesics from $m$ to $a,b,c$ and those from $0$ to $A,B,C$ resp.
We claim that it is not possible to have $\|A-B\|> d(a,b), \|B-C\|> d(b,c),\|A-C\|> d(a,c)$ simultaneously. \\
Let $\bigtriangleup(\bar{a},\bar{b},\bar{c})$ be a comparison triangle in $\R^2$ for $\bigtriangleup(abc)$ and $\bar{m}$ the unique minimizing point in the definition of $\bar{\rho}:=\rho(\bar{a},\bar{b},\bar{c})$. If $\|A-B\|> d(a,b), \|B-C\|> d(b,c)$ and $\|A-C\|> d(a,c)$, then $
 \angle_{\bar{m}}(\bar{a},\bar{b})<\angle_0(A,B)=\angle_{0_m}(V_a,V_b)$ and so on, because $\|\bar{a}-\bar{m}\|=\bar{\rho}r_a\geq\rho r_a=\|A\|$ etc, where $r_a+r_b=d(a,b)$, $r_a+r_c=d(a,c)$ and $r_b+r_c=d(b,c)$. Therefore, the sum of  the angles at $\bar{m}$ is strictly smaller than the sum of the  angles at $0_m$, which itself is not greater than $2\pi$, and this is a contradiction.\\
 Thus, at least one of the side lengths in the triangle $(a,b,c)$ is $\ge$ the corresponding one in the Euclidean triangle with the same angles at the center. 
This holds for arbitrarily small triangles, and therefore, we get a corresponding comparison property between the tangent plane in $T_mN$ spanned by $V_a,V_b$ and $V_c$.
The infinitesimal expansion property characterizes the sectional curvature of that plane and is independent of the directions in that plane. That is, if it holds for two tangent vectors with a nonzero angle between them, then it holds for any pair of vectors in that plane. Thus, we have verified the infinitesimal expansion property characteristic of nonpositive sectional curvature. \qed
\end{pf}

\section{Tripod spaces}\label{tripods}
According to Definition \ref{tripod space} and its subsequent discussion, the metric space $(X,d)$ is a tripod space if and only if $\rho(x_1,x_2,x_3)=1$ and the infimum is attained by some points, each of which is the center of a tripod with legs $x_1, x_2$ and $x_3$, for any triple of points $\{x_1,x_2,x_3\}$ in $X$.
\begin{lemma}
  In a tripod space, the shortest geodesics from $x_i$ to $m$ and from $x_j$ to $m$ form a geodesic from $x_i$ to $x_j$, for $1\leq i<j\leq 3$. 
\end{lemma}
\begin{pf}
  This follows from $d(x_i,m)+d(x_j,m)=d(x_i,x_j)$.\qed
\end{pf}
\begin{theo}
Let $\{(X_i,d_i)\}_{i\in I}$ be a family of tripod spaces and for each $i\in I$ fix a point $z_i\in X_i$. The $\mathit{l}_1$ product of this family with the chosen fixed points, i.e. the set of all $(x_i)_{i\in I}$ with $x_i\in X_i$ and $\sum_{i\in I}d_i(x_i,z_i)<\infty$ endowed with the metric $d_{1}((x_i),(y_i)):=\sum_{i\in I}d_i(x_i,y_i)$, is a tripod space.
\end{theo}
\begin{pf}
Let $X:=\{(x_i)_{i\in I}:\: x_i\in X_i, \sum_{i\in I}d_i(x_i,z_i)<\infty\}$ and $x=(x_i)$, $y=(y_i)$ and $w=(w_i)$ be three points in $X$. Moreover, let $\{r_1, r_2,r_3\}$ be the unique solution of the system $r_1+r_2=d_1(x,y)$, $r_1+r_3=d_1(x,w)$ and $r_2+r_3=d(y,w)$ obtained by the Gromov products \rf{Gromov product}. For each $i\in I$, considering the Gromov products 
\begin{align*}
 \alpha_i&=\frac{1}{2}(d(x_i,y_i)+d(x_i,w_i)-d(y_i,w_i)),&\\
 \beta_i&=\frac{1}{2}(d(x_i,y_i)+d(y_i,w_i)-d(x_i,w_i)),&\\
\gamma_i&=\frac{1}{2}(d(x_i,w_i)+d(y_i,w_i)-d(x_i,y_i)),&
\end{align*} 
there exists at least one point $q_i\in X_i$ such that $d(x_i,q_i)=\alpha_i$, $d(q_i,y_i)=\beta_i$ and $d(w_i,q_i)=\gamma_i$, since $X_i$ is a tripod space for $i\in I$. Therefore, $d_1((x_i),(q_i))=r_1$, $d_1((y_i),(q_i))=r_2$ and $d_1((w_i),(q_i))=r_3$, i.e. $q=(q_i)_{i\in I}$ is the center of a tripod with legs $x$, $y$ and $w$. \qed
\end{pf}
\begin{coro}
For an arbitrary index set $I$, taking $\R$ as $X_i$ and $0$ as the fixed point for all $i\in I$, then the corresponding $\mathit{l}_{1}$ product, i.e. the Banach space $\mathit{l}_{1}(I)$, is a tripod space since $\R$ is a tripod space. The space $\mathit{l}_{1}(I)$ is $(\R^n,d_1)$, when $I$ is a finite set with $n$ elements.
\end{coro}
\begin{theo}\label{tripod and hyperbolicity}
Let $(X,d)$ be a complete metric space that satisfies the tripod condition \rf{0}. Then { $X$ is a geodesic space and }for any three points $x_1,x_2,x_3\in X$ and each 
$\epsilon > 0$ there is an $\epsilon$-slim triangle with vertices $x_1$, $x_2$ and $x_3$.
\end{theo}
\begin{pf}
According to \rf{rho}, for each positive $\epsilon$, there exists a point $\alpha\in X$ such that $d(\alpha,x_i)\leq r_i+\frac{\epsilon}{6}$. Put $y_1=\alpha$ and by induction define $y_{k+1}$ as a point in $X$ for which 
\begin{align*}
\begin{cases} d(y_{k+1},y_k)\leq \frac{\epsilon/3}{2^k}+\frac{\epsilon/3}{2^{k+1}} &\\
d(y_{k+1},x_i)\leq r_i+\frac{\epsilon/3}{2^{k+1}}& \text{for i=1,2}.
\end{cases}
\end{align*}  
The sequence $\{y_k\}$ is Cauchy and therefore there exists $y\in X$ such that $y_k$ approaches $y$ as $k\rightarrow\infty$. Moreover, we have $d(x_i,y)=r_i$ for $i=1,2$, which in particular shows that every tripod space is geodesic. Hence there is a geodesic $c_{12}$ connecting $x_1$ and $x_2$ passing through $y$. We can build a geodesic by attaching any geodesic connecting $x_1$ to $y$ to any geodesic connecting $y$ to $x_2$.\\
By the same process, starting from $\alpha$, we can find a point $z$( resp. $w$) in $X$ and a geodesic $c_{13}$( resp. $c_{23}$) connecting $x_1$ and $x_3$( resp. $x_2$ and $x_3$) passing through $z$( resp. w). The insize of the triangle with sides $c_{12}$, $c_{13}$ and $c_{23}$, i.e. the diameter of $\{y,z,w\}$, is  $\le \epsilon$ by the triangle inequality and therefore this triangle is $\epsilon$-slim.  \qed
\end{pf}
It is worth mentioning that we can weaken the assumption of the theorem above to $\rho(x_1,x_2,x_3)=1$ for each triple $\{x_1,x_2,x_3\}\subset X$ and still get the same result.  
\begin{defi}
Let $(X,d)$ be a tripod space, and let $c_0, c_1:[0,1]\longrightarrow X$ be two shortest geodesics of length $a$ between the same endpoints $c_0(0)=c_1(0)=:x_0$ and $c_0(1)=c_1(1)=:x_1$. A {\em geodesic homotopy} between $c_0$ and $c_1$ is a continuous family $c:[0,1]\times [0,1]\to X$ with $c(t,0)=c_0(t), c(t,1)=c_1(t)$ and with the property that the curves $c(.,s)$ are shortest geodesics from $x_0$ to $x_1$, for all $s\in [0,1]$.  
\end{defi}
In particular, all the geodesics $c(.,s)$ then have the same length.
\begin{theo}\label{hom}
 Let $(X,d)$ be a tripod space. Then any two shortest geodesics $c_0,c_1$ connecting two fixed points are contained in a geodesic homotopy.  
\end{theo}
\begin{pf}
Let $x_0=c_0(0)=c_1(0)$, $x_1=c_0(1)=c_1(1)$ and $a:=d(x_0,x_1)$, which is equal to the lengths of $c_0$ and $c_1$. 
For any $t\in [0,1]$, we then have
\ba
\label{1}
d(c_0(t),c_1(t))\le& \begin{cases} 2at &\text{ for } t\le \frac{1}{2}\\
2a(1-t) &\text{ for } t\ge \frac{1}{2}
\end{cases}\\
\label{2}
\le& a
\ea
and in general, this is smaller than $a$, unless $t=1/2$ and $c_0$ and $c_1$ together form a closed geodesic.

We may assume without loss of generality  that $c_0$ and $c_1$ have only their endpoints $c_0(0)=c_1(0)$ and $c_0(1)=c_1(1)$ in common, but are otherwise disjoint. 
We then consider a tripod with vertices $c_0(2/3),c_1(2/3)$ and $c_0(0)=c_1(0)$. Since by the shortest property and  \rf{1}, all distances are $\le 2a/3$, the tripod cannot coincide with the configuration $c_0,c_1$. Note that he center of the tripod might lie on $c_0$ or $c_1$, and one of the two geodesics forming the tripod might coincide with a portion of $c_0$ or $c_1$, but this cannot happen for both of them simultaneously, as the center of the tripod cannot be $c_0(0)=c_1(0)$, because the geodesic connecting $c_0(2/3)$ and $c_1(2/3)$ through that center has length at most $2a/3$. Therefore, we have a geodesic arc of length $\frac{2}{3}a$ not coinciding with a portion of $c_0$ or $c_1$ from either $c_0(2/3)$ or $c_1(2/3)$ to $c_0(0)=c_1(0)$. That is, we have two shortest geodesics  from $x_0$ to $x_1$ going through the tripod to $c_0(2/3)$ or $c_1(2/3)$ and thence along $c_0$ or $c_1$, resp., to $x_1$. We call them $c_{1/3}$ and $c_{2/3}$. At least one of them does not coincide with $c_0$ or $c_1$. The paths $c_{1/3}$ and $c_{2/3}$ then also have length $a$ and for $t\geq 2/3$ both $d(c_0(t),c_{1/3}(t))$ and $d(c_1(t),c_{2/3}(t))$ are zero. Moreover, the restriction of $c_0$ and $c_{1/3}$ to $0\leq t\leq 2/3$ form two geodesics between $x_0$ and $c_0(2/3)$ of lengths $\frac{2}{3}a$ and so with the same argument as in \ref{1} we have 
\ba
\label{3}
d(c_0(t),c_{1/3}(t)) \le& \begin{cases}  2at &\text{ for } t\le \frac{1}{3}\\
 2a(2/3-t) & \frac{1}{3}\le t\le \frac{2}{3}.
\end{cases}\\
\le & \frac{2}{3}a
\ea
The same is correct for $d(c_1(t),c_{2/3}(t))$, and similarly also
\bel{4}
d(c_{1/3}(t),c_{2/3}(t)) \le \frac{2}{3}a.
\qe
In fact,  for $t\geq 2/3$, the two geodesics $c_{1/3}$ and $c_0$ coincide, and so do $c_{2/3}$ and $c_{1}$ and therefore  $d(c_{1/3}(t),c_{2/3}(t))\le 2a/3$ in this interval. The other part of either $c_{1/3}$ and $c_{2/3}$ consists of two segments, one from $x_0$ to $m$ and the other from $m$ to $c_0(2/3)$ and $c_1(2/3)$ respectively. The latter coincide, and since the length of the geodesic connecting $c_0(2/3)$ and $c_1(2/3)$ via $m$ has length at most $\frac{2}{3}a$, \rf{4} follows.\\
We can then repeat this procedure, replacing 
the geodesic $c_0$ (resp. $c_1$)  by the part of $c_0$ (resp. $c_{1/3}$) that connects $x_0$ to $c_0(2/3)$ and attaching the resulting geodesics to $c_0|_{[2/3,1]}$ to find two geodesics between $c_{0}$ and $c_{1/3}$ with distance $\leq\frac{2}{3}\frac{2}{3}a$ from each other and from the geodesics they were constructed from. Similarly, we can find two geodesics between   $c_1$ by $c_{2/3}$ ( resp. $c_{1/3}$ and $c_{2/3}$) with the same distance bound. 
Iterating, we find families of geodesics between $c_0$ and $c_1$, connecting $x_0$ and $x_1$ with arbitrarily small distances between them. Passing to the limit, we obtain a continuous family of shortest geodesics between $x_0$ and $x_1$, interpolating between $c_0$ and $c_1$.\qed
\end{pf}

\begin{coro}
  The space of shortest geodesics between two points in a tripod space is contractible.
\end{coro}
\section{Conclusions}
The well known \v{C}ech construction assigns  to any cover $\mathcal{U}=(U_i)_{i \in I}$ of a topological space $X$ a simplicial complex $\Sigma(\mathcal{U})$ with vertex set $I$ and a simplex $\sigma_J$ whenever $\bigcap_{j\in J} U_j\neq \emptyset$ for $J\subset I$. When all those intersections are contractible, the homology of $\Sigma(\mathcal{U})$ equals that of $X$ (under some rather general topological conditions on $X$). When $(X,d)$ is a metric space, it is natural to use covers by distance balls. Now, when $(X,d)$ is a {\it hyperconvex} metric space, and if we use a cover $\mathcal{U}$ by distance balls, then whenever 
\bel{cc1}
\bigcap_{j\in J\backslash \{j_0\}} U_j\neq \emptyset \text{ for every }j_0\in J,
\qe
then also
\bel{cc2}
\bigcap_{j\in J} U_j\neq \emptyset,
\qe
that is, whenever $\Sigma(\mathcal{U})$ contains all the boundary facets of some simplex, it also contains that simplex itself. That is, it has no holes, and all its homology groups (except $H_0$, of course) vanish. Therefore, hyperconvex spaces have trivial \v{C}ech homology.\\
When we relax the hyperconvexity to $\lambda$-hyperconvexity for $\lambda >1$ (or $0$-hyperbolicity to $\delta$-hyperbolicity for $\delta>0$), then nontrivial homology groups may arise. From that perspective, hyperconvex spaces are the simplest model spaces, and homology can be seen as a topological measure for the deviation from such a model. \\
Homology groups lead to Betti numbers as integer invariants of topological spaces. Geometry, in contrast, can provide more refined real valued invariants. And after \cite{Riemann:1868}, the fundamental geometric invariants are curvatures. As we have seen in the preceding, the essential geometric content of curvature can be extracted for general metric spaces. From our discussion, it emerges that the basic class of model spaces for curvature  is given by the tripod spaces, a special class of hyperconvex spaces. (Since curvature is a 3-point condition, instead of the general hyperconvexity condition, for our geometric purposes, it is natural to work with a corresponding 3-point version, and that leads to the definition of a tripod space.) From that perspective, the geometric content of curvature in the abstract setting considered here is the deviation from the tripod condition. Euclidean spaces, in contrast, are relegated only to a subsidiary role, based on a normalization of curvature that assigns the value $0$ to them. Of course, in our presentation, we also use Euclidean spaces as comparison spaces, but mainly for the purpose of comparing our curvature concept with traditional ones in Section \ref{upper}. \\
Considering Euclidean spaces as model spaces is traditionally  justified by the fact that spaces whose universal cover has  synthetic curvature $\le 0$ in the sense of Alexandrov are homotopically trivial in the sense that their higher homotopy groups vanish. In technical terms, they are $K(\pi,1)$ spaces, with $\pi$ standing for the first homotopy group. The perspective developed here, however, is a homological and not a homotopical one, and therefore, our natural comparison spaces are tripods. We start their investigation in Section \ref{tripods}. A more systematic investigation of their properties should be of interest. \\
In order to get stronger topological properties, like those of hyperconvex spaces, which are homologically trivial, we might need conditions involving collections of more than three points.\\ 
In fact, If $X$ is a tipod Banach space on which every collection of four closed balls $\{B(x_i,r_i)\}_{i=1}^4$ with non-empty pairwise intersection has a non-void intersection, then every sequence of closed balls with non-empty pairwise intersection has also a non-trivial intersection, c.f. \cite{Lindenstrauss:1962}

\end{document}